\documentclass[a4paper,10pt]{article}
\title{Uniform Distribution of Sequences and its interplay with Functional Analysis}

\author{S.K.Mercourakis  \and G.Vassiliadis}

\date{}

\usepackage{amsmath}
\usepackage{amsfonts}
\usepackage{amssymb}
\usepackage{amscd}
\usepackage{amsthm}
\usepackage{makeidx}
\usepackage{euscript}
\usepackage{graphicx}

\theoremstyle{plain}
\newtheorem{theo}{Theorem}

\newtheorem{lemm}{Lemma}
\newtheorem{prop}{Proposition}
\newtheorem{cor}{Corollary}
\newtheorem{rem}{Remark}
\theoremstyle{definition}
\newtheorem{defn}{Definition}

\DeclareMathOperator\supp{supp}
\DeclareMathOperator\co{co}
\DeclareMathOperator\ex{ex}
\DeclareMathOperator\diam{diam}
\DeclareMathOperator\R{Re}
\DeclareMathOperator\I{Im}

\begin{document}
\maketitle

\begin{abstract}
\normalsize
In this paper we apply ideas from the theory of Uniform Distribution of sequences to Functional Analysis
and then drawing inspiration from the consequent results, we study concepts and results in Uniform Distribution itself.
So let $E$ be a Banach space. Then we prove:\\
(a) If $F$ is a bounded subset of $E$ and $x \in \overline{\co}(F)$ (= the closed convex hull of $F$), then there is
a sequence $(x_n) \subseteq F$ which is Ces\`{a}ro summable to $x$.\\
(b) If $E$ is separable, $F \subseteq E^*$ bounded and $f \in \overline{\co}^{w^*}(F)$, then there is a sequence
$(f_n) \subseteq F$ whose sequence of arithmetic means $\frac{f_1+\dots+f_N}{N}$, $N \ge 1$ weak$^*$-converges to $f$.

By the aid of the Krein-Milman theorem, both (a) and (b) have interesting implications for closed, convex and bounded
subsets $\Omega$ of $E$ such that $\Omega=\overline{\co}(\ex \Omega)$ and for weak$^*$ compact and convex subsets of $E^*$.
Of particular interest is the case when $\Omega=B_{C(K)^*}$, where $K$ is a compact metric space.

By further expanding the previous ideas and results, we are able to generalize a classical theorem of Uniform Distribution
which is valid for increasing functions $\varphi:I=[0,1] \rightarrow \mathbb{R}$ with $\varphi(0)=0$ and $\varphi(1)=1$, for functions $\varphi$
of bounded variation on $I$ with $\varphi(0)=0$ and total variation  $V_0^1 \varphi=1$.

\footnote{\noindent keywords: {uniformly distributed sequence \and Ces\`{a}ro summable sequence \and function of bounded variation.}
{Primary  46B09 \and 11K06; Secondary  40C05 \and 60B10.}}

\end{abstract}

\large

\section*{Introduction}

Our aim in this paper is twofold. We first study consequences of ideas coming from the theory of Uniform Distribution of sequences \cite{KN} 
in Functional Analysis (section 1) and then we investigate concepts and results of the theory of Uniform Distribution itself,
setting them in a more general framework (section 2).

In the first section we generalize and improve an important result of Niederreiter \cite{N}, which we state as Theorem 1.
The translation of this theorem into the language of Functional Analysis is Theorem 2, which is the main result of this section.
The first assertion of this theorem is strongly related to a celebrated lemma due to Maurey (see Lemma D of \cite{BPST})
and roughly says that given a point $x$ in the convex hull $\co(F)$ of a bounded subset of some Banach space, then for every $N \in \mathbb{N}$,
$x$ can be approximated (an estimation of the approximation error is also given) by the arithmetic mean of $N$ points of $F$.
The second assertion says that, if $x$ belongs to $\overline{\co}(F)$, then there is a sequence of points of $F$ which is Ces\`{a}ro
summable to $x$. A number of easy consequences of Theorem 2 for a Banach space $E$ are the following:\\
A) If $(y_n) \subseteq E$ is any weakly null sequence, then there is a function $\varphi:\mathbb{N} \rightarrow \mathbb{N}$, such that
the sequence $x_n=y_{\varphi(n)},\; n \ge 1$ is Ces\`{a}ro summable to zero (Proposition 2).\\
B) If $\Omega$ is a closed convex bounded subset of $E$ equal to the closed convex hull of its extreme points $\ex \Omega$, then
for every $x \in \Omega$ there is a sequence $(x_n) \subseteq \ex \Omega$ which is Ces\`{a}ro summable to $x$ (Proposition 3). This 
result has, by the aid of the Krein-Milman theorem, obvious implications for the unit ball of a Banach space which is reflexive or 
of the form $C(K)$ (=the space of real continuous functions on $K$), where $K$ is any compact totally disconnected space (Corollaries 1 and 2).\\
C) A result analogous to Theorem 2(2) (with  analogous proof) is valid for the dual $(E^*,weak^*)$ of a separable Banach space $E$.
When $F$ is a bounded subset of $E^*$ and $f$ belongs to $\overline{\co}^{w^*}(F)$, there is a sequence $(f_n) \subseteq F$
whose arithmetic means weak$^*$-converge to $f$ (Proposition 4).

This result (again using the Krein-Milman theorem) has obvious implications for a weak$^*$-compact and convex subset $\Omega$ of $E^*$,
which generalize classical results (see Proposition 5, Corollary 3, Theorem 3 and 4). So Theorem 3 is the well known result stating that every 
probability measure $\mu$ on a compact metric space $K$ admits a uniformly distributed (u.d.) sequence, but Theorem 4 says something that seems to be new:
For every signed measure $\mu \in \Omega=B_{M(K)}$, there is a sequence $(x_n) \subseteq K$ and a sequence of signs $(\varepsilon_n)
\subseteq \{\pm 1\}$, so that the sequence $\mu_N=\frac{\varepsilon_1 \delta_{x_1}+ \dots +\varepsilon_N \delta_{x_N}}{N},\; N \ge 1$,
weak$^*$-converges to $\mu$.

The last result is our motivation for the second section of this paper. Drawing inspiration from Theorem 4, we extend the classical
concept of uniformly distributed sequence defined for probability measures $\mu$ on a compact space $K$ (see Definition 1.1 of \cite{KN}),
to every (real) signed measure $\mu \in M(K)$ with $\|\mu\|=1$. Thus we will say that a sequence $(x_n) \subseteq K$ is $\mu$-u.d. iff
$(x_n)$ is $|\mu|$-u.d. (in the classical sense) and if also there is a sequence of signs $(\varepsilon_n) \subseteq \{\pm 1\}$,
so that the sequence $\mu_N=\frac{\varepsilon_1 \delta_{x_1}+ \dots +\varepsilon_N \delta_{x_N}}{N},\; N \ge 1$
weak$^*$-converges to $\mu$ (Definition 1).

Then (generalizing Theorem 3) we prove Theorem 5, which states that given a compact metric space $K$ and $\mu \in M(K)$ with $\|\mu\|=1$,
then $\mu$ admits a u.d. sequence $(x_n) \subseteq K$ (in the sense of the aforementioned definition). So if $f \in C(K)$, we have
\[\lim_{N \to \infty} \frac{f(x_1)+\dots+f(x_N)}{N}=\int_K f d |\mu| \;\, \textrm{and} \;\, 
\lim_{N \to \infty} \frac{\varepsilon_1 f(x_1)+\dots+ \varepsilon_N f(x_N)}{N}=\int_K f d \mu.\]
A further generalization of the last theorem is Theorem 6, which says that both of the above equalities are also valid for
$\mu$-Riemann integrable functions. Now let $K$ be a compact interval of the real line, say for simplicity $K=I=[0,1]$.
Taking into account the standard identification of signed Borel measures on $I$ with (proper) functions of bounded variation (BV) on $I$,
Theorem 6 yields Theorem 7: 
Let $\varphi:I \rightarrow \mathbb{R}$ be a BV function with $\varphi(0)=0$, $V_0^1 \varphi=1$ and $\varphi$ is right continuous on $I$.
Then there are sequences $(x_n) \subseteq I$ and $(\varepsilon_n) \subseteq \{\pm 1\}$, such that for every point of continuity $x$ of $\varphi$ we have
\[ \lim_{N \to \infty} \frac{1}{N} \sum_{k=1}^N \chi_{[0,x)}(x_k)=\upsilon(x) \;\, \textrm{and} \;\,  \lim_{N \to \infty} \frac{1}{N} 
\sum_{k=1}^N \varepsilon_k \chi_{[0,x)}(x_k)=\varphi(x),\]
where $\upsilon$ is the function of total variation of $\varphi$ on $I$ and $\upsilon(1)=V_0^1 \varphi$. The last theorem partially generalizes
a classical result from \cite{KN} (Theorem 8 in our treatment) which says that the equalities of Theorem 7 are valid for every $x \in I$,
provided that $\varphi$ is increasing with $\varphi(0)=0$ and $\varphi(1)=1$ (of course, since $\varphi$ is increasing, we have $\upsilon=\varphi$
and $\varepsilon_k=1$ for all $k \ge 1$).

The rest of this section is devoted to the proof of Theorem 9, that is, of the fact that Theorem 7 holds true for every BV function $\varphi$
on $I$ with $\varphi(0)=0$, $V_0^1 \varphi=1$ and for each point $x \in I$. This result is a common generalization of Theorems 7 and 8 and is the main result of the second section. The proof of Theorem 9 is rather elaborate and is presented in several steps (Lemmas 5, 6, 7 etc.).
We also note that the notion of discrepancy of a sequence in $I$ is crucial in the proof of Theorem 9.

\section*{Preliminaries}
If $E$ is any Banach space, then $B_E$ denotes its closed unit ball. A subset $L$ of $E$ is said to be total in $E$, if its linear span 
$\langle L \rangle$ is dense in $E$. A sequence $(x_n) \subseteq E$ is said to be Ces\`{a}ro summable, if the corresponding
sequence $\frac{x_1+\dots+x_n}{n},\; n \ge 1$ of arithmetic means of $(x_n)$ converges in norm. Let $A \subseteq E$, then $\co(A)$ is
the convex hull of $A$ and $\overline{\co}(A)$ the norm closure of $\co(A)$, which by a classical theorem of Mazur coincides with the weak
closure of $\co(A)$. If $A \subseteq E^*$, then $\overline{\co}^{w^*}(A)$ denotes the closure of $\co(A)$ in the weak$^*$ topology of $E^*$.
Let $x \in \co(A)$ with $x=\sum_{k=1}^n \lambda_k x_k$, where $x_1,\dots,x_n$ are distinct points of $A$, $\lambda_k>0$ for $k=1,2,\dots,n$
and $\sum_{k=1}^n \lambda_k=1$, then we set $\supp x=\{x_1,\dots,x_n\}$.

Let $K$ be a compact Hausdorff space, then $C(K)$ is the Banach space with sup-norm (denoted by $\|\cdot\|_\infty$) of all continuous real valued functions on $K$. 
The dual $C(K)^*$ 
of $C(K)$ is isometrically identified via the classical Riesz representation theorem with the space $M(K)$ of all finite regular
signed Borel measures on $K$, with norm $\|\mu\|=|\mu|(K)$. By $M^+(K)$ (resp. $P(K)$) we denote the positive (resp. probability)
measures on $K$. When $\mu \in M^+(K)$, a bounded function $f:K \rightarrow \mathbb{R}$ is said to be $\mu$-Riemann integrable, if the 
set of discontinuity points of $f$ has $\mu$-measure zero. It is an easy consequence of Lusin's theorem that each $\mu$-Riemann integrable
function is $\mu$-measurable and hence $\mu$-integrable.

Let $X$ be a (nonempty) set; then $|X|$ denotes the cardinality of $X$ and $\ell_\infty(X)$ the Banach space (with sup-norm) of all bounded
real valued functions on $X$. It is well known that $\ell_\infty(X)$ is linearly isometric to the space $C(\beta X)$, where $\beta X$ is the 
Stone-\v{C}ech compactification of the discrete set $X$. If $x \in X$, then $\delta_x$ denotes the point mass at $x$, that is, the Dirac measure
$\delta_x:\ell_\infty(X) \rightarrow \mathbb{R}$ such that $\delta_x(f)=f(x)$, for $f \in \ell_\infty(X)$. We denote by $\mathcal{F}(X)$
the set of probability measures of finite support on $X$, thus $\mathcal{F}(X)=\co(\{\delta_x: x \in X\})$; if $\mu \in \mathcal{F}(X)$,
$x_1,\dots,x_n$ are distinct points of $X$ such that $\mu(\{x\})=\lambda_k>0$, for $k=1,2,\dots,n$ and $\sum_{k=1}^n \lambda_k=1$,
then $\mu=\sum_{k=1}^n \lambda_k \delta_{x_k}$ and thus $\supp \mu=\{x_1,\dots,x_n\}$. Also, when $A \subseteq X$ we denote by $\chi_A$ the
characteristic function of $A$. We note that if $X$ is compact Hausdorff, then $P(X)$
is weak$^*$ compact and convex subset of $M(X)=C(X)^*$ and the set of its extreme points $\ex P(X)$ coincides with the set of Dirac
measures on $X$; therefore by Krein-Milman's theorem $P(X)=\overline{\co}^{w^*}(\{\delta_x:x \in X\})=\overline{\mathcal{F}(X)}^{w^*}$.

Let $X$ be compact Hausdorff and $\mu \in P(X)$. A sequence $(x_n) \subseteq X$ is called $\mu$-uniformly distributed (shortly $\mu$-u.d.)
in $X$ if
\[ \lim_{N \to \infty} \frac{1}{N} \sum_{k=1}^N f(x_k) =\int_X fd \mu  \;\; \textrm{for all} \;\; f \in C(X) \]
(equivalently  if weak$^*-\lim_{N \to \infty} \frac{\delta_{x_1}+\dots+\delta_{x_N}}{N}=  \mu$.)

Whilst most of our results remain valid in the complex case, we assume for simplicity that all Banach spaces (and functions) are real
and in certain cases we indicate what happens in the complex case.

\section{Functional Analytic consequences of a result of Niederreiter}

We begin by generalizing and improving an important result of Niederreiter,
essentially following the proof of the original result (see Theorem 1 of \cite{N}).

\begin{theo}
Let $X$ be a nonempty set, $L$ a subset of the closed unit ball $B$ of $\ell_\infty(X)$
and $(\mu_j) \subseteq \mathcal{F}(X)$.
Assume that the sequence $(\mu_j)$ converges pointwise on $L$, that is, there exists a function
$\mu:L \rightarrow \mathbb{R}$ such that 
\[\mu_j(f) \underset{j \to \infty}{\longrightarrow} \mu(f) \;\;\; \forall f \in L.\]
Then there is a sequence $\omega=(x_n) \subseteq \bigcup_{j=1}^\infty \supp \mu_j$ such that
\begin{enumerate}
\item The sequence $\upsilon_N=\frac{\delta_{x_1}+\dots+\delta_{x_N}}{N} \underset{N \to \infty}{\longrightarrow} \mu$
pointwise on $L$.
\item Moreover, if the sequence $(\mu_j)$ converges to $\mu$ uniformly on $L$, then $(\upsilon_N)$
converges to $\mu$ uniformly on $L$.
\end{enumerate}
\end{theo}

The main tool for proving the above theorem is the following lemma (see Lemma 1 of \cite{N}).

\begin{lemm}
Let $\mu \in \mathcal{F}(X)$; then there exists a positive constant $C(\mu)$ and a sequence $\omega=(y_n)$
in $X$, such that

\begin{equation}
\left| \frac{1}{N}\sum_{k=1}^N \chi_M(y_k)-\mu(M) \right| \le \frac{C(\mu)}{N} 
\end{equation}

for all $N \in \mathbb{N}$ and for all subsets $M \subseteq X$. In particular
\[C(\mu)=(m-1)\left[\frac{m}{2}\right] \]
will do, where $m=| \supp \mu |$.

\end{lemm}

It is necessary for our purposes to prove that a modification of the above lemma holds,
not only for characteristic functions, but also for every bounded function $f:X \rightarrow
\mathbb{R}$. We recall that the set of extreme points $\ex B$ of the unit ball $B$ of 
$\ell_\infty(X)$ consists of all functions $f:X \rightarrow
\mathbb{R}$ such that $|f(x)|=1$, for all $x \in X$ and the well known fact that 
$B=\overline{\co}(\ex B)$.

\begin{prop}
Let $f:X \rightarrow \mathbb{R}$ be any bounded function. Then with positive constant $C(\mu)$ and 
sequence $\omega=(y_n)$ in $X$ of Lemma 1, inequality (1)
holds in the following modified form.\\
(a) If $f \in \co(\ex B)$ then we have

\begin{equation}
\left| \frac{1}{N}\sum_{k=1}^N f(y_k)-\int_X fd \mu \right| \le 2 \frac{C(\mu)}{N} 
\end{equation}
for all $N \ge 1$.\\
(b) If $f$ is any bounded function, then we have

\begin{equation}
\left| \frac{1}{N}\sum_{k=1}^N f(y_k)-\int_X fd \mu \right| \le  \frac{2 \|f\|_\infty}
{N} (1+C(\mu)) 
\end{equation}
for all $N \ge 1$.

\end{prop}

\begin{proof}  
(a) Assume first that $f \in \ex B$. Set $V=\{x \in X:f(x)=1\}$, then the complement
of $V$ is the set $V^c=\{x \in X:f(x)=-1\}$. Therefore $f=\chi_V-\chi_{V^c}$. So we get 
for $N \in \mathbb{N}$ that
$\sum_{k=1}^N f(y_k)=\sum_{k=1}^N \chi_V(y_k)-\sum_{k=1}^N \chi_{V^c}(y_k)$ and 
$\int_X f d\mu=\mu(V)-\mu(V^c)$.
Now from Lemma 1 we have,
\[\left| \frac{1}{N}\sum_{k=1}^N f(y_k)-\int_X fd \mu \right| \le
\left| \frac{1}{N}\sum_{k=1}^N \chi_V(y_k)-\mu(V) \right| +
\left| \frac{1}{N}\sum_{k=1}^N \chi_{V^c}(y_k)-\mu(V^c) \right|\]
\[ \le \frac{C(\mu)}{N}+\frac{C(\mu)}{N}=2\frac{C(\mu)}{N}. \]
It now follows easily from the last inequality that (2) remains valid, for all
$f \in \co(\ex B))$.

(b) It is clear that it suffices to prove (3) for $f \in B$. Since $B=\overline{\co}(\ex B)$,
there is a sequence $(f_n) \subseteq \co(\ex B)$ such that $f_n \rightarrow f$ uniformly on $X$.
Given $N \in \mathbb{N}$, consider $n_0 \in \mathbb{N}$ such that

\begin{equation}
\|f-f_{n_0}\|_\infty<\frac{1}{N}.
\end{equation} 

Then from assertion (a) and (4) we get that
\begin{equation*}
\left| \frac{1}{N}\sum_{k=1}^N f(y_k)-\int_X fd \mu \right| 
\end{equation*}
\begin{equation*}
=\left| \frac{1}{N}\sum_{k=1}^N 
(f(y_k)-f_{n_0}(y_k))+ \left(\int_X f_{n_0} d\mu-\int_X fd \mu \right)+ 
\left(\frac{1}{N}\sum_{k=1}^N f_{n_0}(y_k)-\int_X f_{n_0} d \mu\right) \right| 
\end{equation*}
\begin{equation*}
 \le \frac{1}{N}\sum_{k=1}^N |f(y_k)-f_{n_0}(y_k)|+ \int_X |f-f_{n_0}| d\mu+ 
\left|\frac{1}{N}\sum_{k=1}^N f_{n_0}(y_k)-\int_X f_{n_0} d \mu\right| 
\end{equation*}
\begin{equation*}
\le \frac{1}{N} \cdot  N \cdot  \frac{1}{N}+\frac{1}{N} \mu (X)+2 \frac{C(\mu)}{N}=
\frac{2}{N} (1+C(\mu)).
\end{equation*}

\end{proof}

\begin{rem}
Let $f=\R f+i \I f$ be any complex function such that\\ $|f(x)|=\sqrt{(\R f(x))^2+(\I f(x))^2} \le 1$
for all $x \in X$. Then it is easy to prove that
\begin{equation}
\left| \frac{1}{N}\sum_{k=1}^N f(y_k)-\int_X fd \mu \right| \le
\frac{2\sqrt{2}}{N} (1+C(\mu)) 
\end{equation}
for all $N \ge 1$.

In particular (5) is valid for any extreme point $f$ of the unit ball $B$ of the complex
Banach space $\ell_\infty(X)$ (recall that the extreme points of $B$ are the functions
of the form $f:X \rightarrow \mathbb{C}$, such that $|f(x)|=1$ for all $x \in X$).
\end{rem}

We now proceed with the proof of Theorem 1. Note that assertion (a) is slightly 
more general than Theorem 1 of \cite{N}; assertion (b) is new.

\begin{proof} \textit{(of Theorem 1)}
(1) Assume that $(\mu_j)$ converges pointwise on $L$ to $\mu$. By Lemma 1
there exist positive constants $C_j=C(\mu_j)$ and sequences $\omega_j=(x_{j,n})_{n \ge 1},\; j \in \mathbb{N}$
such that relation (1) of Lemma 1 holds. For each $j \in \mathbb{N}$, choose a positive 
integer $r_j$ such that $r_j \ge \max\{j^2,j(C_1+\cdots+C_{j+1})\}$. Put $r_0=0$; we define a sequence
$\omega=(x_n)$ as follows. Every positive integer $n$ has a unique representation of the form
$n=r_0+r_1+\cdots+r_{j-1}+s$ with $j \ge 1$ and $0<s \le r_j$; we set $x_n=x_{j,s}$. Take an integer
$N>r_1$; $N$ can be written in the form $N=r_1+\cdots+r_k+s$ with $0<s \le r_{k+1}$.
For any function $f \in L$ we get
\[\sum_{n=1}^N f(x_n)=\sum_{j=1}^k \left(\sum_{\lambda=1}^{r_j} f(x_{j,\lambda}) \right)
+\sum_{\lambda=1}^s f(x_{k+1,\lambda}) .\]
Therefore 
\begin{equation*}
|\upsilon_N(f)-\mu(f)|=\left| \frac{1}{N} \sum_{n=1}^N f(x_n)-\mu(f) \right|
\end{equation*}
\begin{equation*}
=|\sum_{j=1}^k \frac{r_j}{N} 
\left( \frac{1}{r_j}\sum_{\lambda=1}^{r_j} f(x_{j,\lambda})-\mu_j(f) \right) + \frac{s}{N} \left( \frac{1}{s}
\sum_{\lambda=1}^s f(x_{k+1,\lambda})-\mu_{k+1}(f) \right)+
\end{equation*}

\begin{equation*}
+\sum_{j=1}^k \frac{r_j}{N} \mu_j(f)+\frac{s}{N}
\mu_{k+1}(f)-\mu(f)|
\end{equation*}
\textrm{(using Proposition 1)} 
\begin{equation*}
\le  \sum_{j=1}^k \frac{r_j}{N} 
\left[ \frac{2}{r_j}(1+C_j)\right] + \frac{s}{N} \cdot \frac{2}{s} (1+C_{k+1})+
\left|\frac{1}{N} \left[\sum_{j=1}^k r_j \mu_j(f)+s \mu_{k+1}(f)\right]-\mu(f) \right|  
\end{equation*}
\begin{equation*}
\le \left(\textrm{by letting}\; K(N,f)= \frac{1}{N} \left[\sum_{j=1}^k r_j \mu_j(f)+s \mu_{k+1}(f)\right]-\mu(f) \right)
\end{equation*}
\begin{equation*}
\le \frac{2}{r_k} \sum_{j=1}^{k+1}(1+C_j)+|K(N,f)|=\frac{2}{r_k}(k+1)+\frac{2}{r_k} \sum_{j=1}^{k+1}C_j+|K(N,f)| 
\end{equation*}
(since $r_k \ge \max\{k^2,k(C_1+\cdots+C_{k+1})\}$)
\begin{equation*}
\le \frac{2(k+1)}{k^2}+\frac{2}{k}+|K(N,f)|.
\end{equation*}
If $N \rightarrow \infty$ then $k \rightarrow \infty$ and the sum of the first two terms tends to zero. In order
to prove that the third term tends to zero, we set for every $N>r_1$
\[A_N=\left(\frac{r_1}{N},\frac{r_2}{N},\cdots,\frac{r_k}{N},\frac{s}{N},0,\cdots \right), \]
where $N=r_1+r_2+\cdots+r_k+s,\; 0<s \le r_{k+1}$. Then $A=(A_N)$ defines an infinite real matrix that is a regular
method of summability. If we set $H_f=(\mu_j(f))_{j \ge 1}$, where $f \in L$, then we have
\[A_N \cdot H_f=\frac{1}{N} \left[\sum_{j=1}^k r_j \mu_j(f)+s \mu_{k+1}(f) \right], \; N \ge 1.\]
Since $\mu_j(f) \underset{j \to \infty}{\longrightarrow} \mu(f)$ for $f \in L$ and $A$ is a regular method of summability,
we get that $|A_N \cdot H_f-\mu(f)|=|K(N,f)| \underset{N \to \infty}{\longrightarrow} 0$ for all $f \in L$ and
we are done.

(2) We assume now that $(\mu_j)$ converges to $\mu$ uniformly on $L$. Since $A$ is a regular method
of summability, we get that
\[A_N \cdot H_f \underset{N \to \infty}{\longrightarrow} \mu(f) \; \textrm{uniformly on} \; L.\]
Therefore given $\varepsilon>0$, there is $N_0=N_0(\varepsilon)$ such that
\[N \ge N_0 \Rightarrow |A_n \cdot H_f-\mu(f)|=|K(N,f)| \le \frac{\varepsilon}{2} \;\; \forall f \in L\]
and of course $2 \left( \frac{k+1}{k^2}+\frac{1}{k} \right) \le \frac{\varepsilon}{2}$, if $N_0$ is
sufficiently large. It then follows from the above that
\[N \ge N_0 \Rightarrow \left| \frac{1}{N} \sum_{n=1}^N f(x_n)-\mu(f) \right| \le \varepsilon \; \;\forall
f \in L,\] 
which means that

\[\upsilon_N=\frac{\delta_{x_1}+\dots+\delta_{x_N}}{N} \underset{N \to \infty}{\longrightarrow} \mu\]
uniformly on $L$.

\end{proof}

\begin{rem}
We notice that using inequality (5) of Remark 1 in the proof of Theorem 1 instead of inequality (3) of Proposition 1,
we can easily prove that Theorem 1 is also valid assuming that $L$ consists of complex functions.
Therefore Theorem 2, which we are going to prove, and everything depending on this theorem is also valid in the
complex case.
\end{rem}

The rest of this section is devoted to some applications of the previous results (Proposition 1 and Theorem 1)
in Banach space theory. We first prove the following

\begin{theo}
Let $E$ be a Banach space, $F$ a bounded subset of $E$ with $F \subseteq B(0,R)$ and $x \in E$.
Then we have:
\begin{enumerate}
\item Assume that $x \in \co(F)$ and let $F_0 \subseteq F$ be any finite set such that $x \in \co(F_0)$.
Then there is a sequence $(x_n) \subseteq F_0$ and a positive constant $C=C(|F_0|)$ such that
\[\left\|\frac{1}{N} \sum_{k=1}^N x_k-x \right\| \le \frac{2R}{N}(1+C) \;\; \textrm{for all} \; N \ge 1;\]
in particular
\[\| \cdot\|-\lim_{N \to \infty} \frac{1}{N} \sum_{k=1}^N x_k =x. \]
\item Assume that $x \in \overline{\co}(F)$. Then there is a sequence $(x_n) \subseteq F$ such that
\[\| \cdot\|-\lim_{N \to \infty} \frac{1}{N} \sum_{k=1}^N x_k =x. \]
\end{enumerate}

\end{theo}

\begin{proof}
Assume without loss of generality that $R=1$, that is $F \subseteq B_E$; otherwise we replace $F$
by $\frac{1}{R} F$ and $x$ by $\frac{1}{R} x$. We set $X=B_E$ and notice that each $f \in E^*$ can be 
identified with a bounded (continuous) function on $X$ through the isometry operator
$T:f \in E^* \mapsto T(f)=f|X \in \ell_\infty(X)$.

(1) Let $x \in \co(F)$, then $x=\alpha_1 y_1+\dots+\alpha_m y_m$, where $y_1,\dots,y_m$ are
distinct points of $F$, $\alpha_k>0$ for $k=1,2,\dots,m$ and $\sum_{k=1}^m \alpha_k=1$. We can
consider $x$ as a finitely supported measure $\mu$ on $X$, by letting 
$\mu=\sum_{k=1}^m \alpha_k \delta_{y_k}$; clearly $\supp \mu=\{y_1,\dots,y_m\}$. Then 
$\mu$ represents $x$, that is, for every $f \in B_{E^*}$
\[\int_X f d\mu=\sum_{k=1}^m \alpha_k f(y_k)=f \left(\sum_{k=1}^m \alpha_k y_k \right)=f(x).\]
It then follows from Proposition 1 (see also Remark 1) that there is $(x_n) \subseteq F_0$, where
$F_0=\supp(\mu)$,  such that for every $f \in \ell_\infty(X)$ with $\|f\| \le 1$ we have
\[\left| \frac{1}{N} \sum_{k=1}^N f(x_k)- \int_X f d\mu \right| \le \frac{2}{N} (1+C(\mu))
\; \textrm{for all} \; N \ge 1.\]
In particular, if $f \in E^*$ with $\|f\| \le 1$, then we have
\[\left| \frac{1}{N} \sum_{k=1}^N f(x_k)- \int_X f d\mu \right|=\left|f \left(\frac{1}{N}
\sum_{k=1}^N x_k -x \right) \right| \le \frac{2}{N} (1+C(\mu))
\; \textrm{for all} \; N \ge 1,\]
which implies that for all $N \ge 1$,
\[ \left\| \frac{1}{N} \sum_{k=1}^N x_k-x \right\| =\sup \left\{
\left|f \left(\frac{1}{N} \sum_{k=1}^N x_k -x \right) \right| :f \in B_{E^*} \right\} 
\le \frac{2}{N} (1+C(\mu)).  \]
We set $C=C(\mu)$, hence $C$ depends on $m=|F_0|$ and obtain the desired result.

(2) Let $x \in \overline{\co}(F)$; then there is a sequence $(\mu_j)$ of convex combinations
of elements of $F$ such that 
\begin{equation}
\| \cdot\|-\lim \mu_j = x. 
\end{equation}

We consider each $\mu_j$ as a finitely supported probability measure on $F \subseteq X=B_E$.
So if we set $L=T(B_{E^*}) \subseteq B_{\ell_\infty(X)}$, then (6) means that
$\mu_j \underset{j \to \infty}{\longrightarrow} x$ uniformly on $L$. It then follows from Theorem 1(2)
that there is a sequence $(x_n) \subseteq \cup_{j=1}^\infty \supp \mu_j \subseteq F$ such that the sequence
of arithmetic means
\[\frac{\delta_{x_1}+\dots+\delta_{x_N}}{N} \underset{N \to \infty}{\longrightarrow} x\;\;
\textrm{uniformly on } L,\]
equivalently $\| \cdot\|-\lim_{N \to \infty} \frac{x_1+\dots+x_N}{N}= x$.
The proof of the theorem is complete.

\end{proof}

Assertion (1) of the above theorem has a strong relationship with a lemma due to Maurey (see Lemma D of \cite{BPST})
which states that

\begin{lemm}
Let $E$ be a Banach space of type $p$ for some $p>1$, $F \subseteq E$ and $x \in \co(F)$. Set 
$q=\frac{p}{p-1}$. Then for every $N \in \mathbb{N}$ there exist $x_1,\dots,x_N \in F$ such that
\[ \left\| \frac{1}{N} \sum_{k=1}^N x_k-x \right\| \le \diam(F) \frac{T_p(E)}{N^{\frac{1}{q}}}.\]
($T_p(E)$ is the type $p$ constant of $E$, see pp. 137-8 of \cite{AK}).
\end{lemm}

In our case the constant $C$ depends on the cardinality of the finite subset of $F$ that supports $x$;
in Maurey's Lemma it depends only on the space $E$, which must be of type $p$. If we assume, as we may,
that $F \subseteq B_E$ (thus $\diam(F) \le 2$) and $N$ is large, then clearly
$\frac{1+C}{N} \le \frac{T_p(E)}{N^{\frac{1}{q}}}$, hence we get Maurey's Lemma. In some way assertion (1)
of Theorem 2 is the "pointwise" version of Maurey's Lemma.

From assertion (2) of Theorem 2 together with Mazur's classical result, stating that the weak and the norm closure
of any convex subset of a Banach space coincide, we obtain the following interesting consequence.

\begin{prop}
Let $E$ be a Banach space and $(y_n) \subseteq E$ be any weakly convergent sequence, so that $y_n \overset{w}{\rightarrow} y$.
Then there is a function $\varphi:\mathbb{N} \rightarrow \mathbb{N}$ such that the sequence $x_n=y_{\varphi(n)}, \;n \ge 1$
is Ces\`{a}ro summable to $y$, that is
$\| \cdot\|-\lim_{N \to \infty} \frac{x_1+\dots+x_N}{N}=y$.
\end{prop}

\begin{proof} 
The set $F=\{y_n:n \ge 1\}$ is bounded. Since $y \in \overline{F}^w$, we get that $y \in \overline{\co}^w(F)$.
By Mazur's theorem we have that $\overline{\co}^w(F)=\overline{\co}(F)$. So $y \in \overline{\co}(F)$
and then assertion (2) of Theorem 2 can be applied.

\end{proof}

\begin{rem}
It is well known that the function $\varphi$ of Proposition 2 cannot in general be chosen strictly increasing (neither 1-1).
In fact, it is possible to find a weakly null sequence $(y_n)$, such that for every subsequence $(y'_n)$ of $(y_n)$
the sequence of arithmetic means $\frac{y'_1+\dots+y'_N}{N}, \; N \ge 1$ is not norm convergent.
The first such example was constructed by J. Schreier (see \cite{M1} and \cite{AMT}). We note in this connection
that in \cite{AMT} is given a complete classification of the complexity of weakly null sequences, by the
aid of a hierarchy of summability methods introduced there.
\end{rem}

Still another immediate but useful consequence of assertion (2) of Theorem 2 is the following.

\begin{prop}
Let $E$ be a Banach space and $\Omega$ be a closed, convex, bounded subset of $E$, such that $\Omega$
is equal to the closed convex hull of its extreme points $\ex \Omega$, that is $\Omega=\overline{\co}(\ex \Omega)$.
Then for every $x \in \Omega$, there is a sequence $(x_n) \subseteq \ex \Omega$ which is Ces\`{a}ro summable to $x$.
\end{prop}

We now present some applications of Proposition 3.

\begin{cor}
Let $\Omega$ be a weakly compact and convex subset of a Banach space $E$ (in particular $\Omega=B_E$ and $E$ is reflexive). 
Then for every $x \in \Omega$ there is a sequence $(x_n)$ of extreme
points of $\Omega$ which is Ces\`{a}ro summable to $x$.
\end{cor}

\begin{proof} 
Since $\Omega$ is weakly compact and convex, by the Krein-Milman theorem we have that
$\Omega=\overline{\co}(\ex \Omega)$. Hence the result is an immediate consequence of Proposition 3.

\end{proof}

The next result concerns Banach spaces of the form $C(K)$, where $K$ is a compact Hausdorff space.
We recall that $K$ is called totally disconnected, if it has a base for its topology consisting of
open and closed (clopen) sets.

\begin{cor}
Let $K$ be a compact Hausdorff space. We assume that either\\
(a) $C(K)$ is the space of continuous complex functions on $K$, or\\
(b) $K$ is totally disconnected and $C(K)$ is the space of continuous real functions on $K$.

Then for every $f \in B=B_{C(K)}$, there exists a sequence $(f_n)$ of extreme points
of $B$ such that 
\[\| \cdot\|_\infty-\lim_{N \to \infty} \frac{f_1+\dots+f_N}{N}= f.\]
\end{cor}

\begin{proof} 
In either case we have that $B=\overline{\co}(\ex B)$ (see Theorems 1.6 and 1.8 of \cite{B}). Hence Proposition 3 can be applied. 

\end{proof}

\begin{rem}
(a) Recall that if $K$ is a compact Hausdorff space, then $f \in \ex B$ iff $|f(x)|=1$, for all $x \in K$ (see Theorem 1.3 of \cite{B}).
If $K$ is in addition totally disconnected, $C(K)$ is the space of real continuous functions on $K$ and $f \in \ex B$, then
the sets $V=\{x \in K:f(x)=1\}$ and $V^c=\{x \in K: f(x)=-1\}$ constitute a partition of $K$ in two clopen sets. Thus the extreme
points of $B$ are completely determined by the clopen nonempty subsets of $K$.\\
(b) We have already used the above remark in the special case of the Banach space $\ell_\infty(X)$ (cf. the proof of Proposition 1).
Indeed $\ell_\infty(X)$ is isometric to $C(\beta X)$, where $\beta X$ is the Stone-\v{C}ech compactification
of the discrete set $X$, which is a compact extremally disconnected space.
\end{rem}

We continue our investigation, applying Theorem 1 to the weak$^*$ topology of the dual $E^*$ of a separable Banach space $E$.
As we shall see, results similar to Theorem 2 (2) and Proposition 3 are valid. Moreover, our approach has interesting applications
for the dual $M(K)=C(K)^*$ of $C(K)$, where $K$ is any compact metric space.

\begin{prop}
Let $E$ be a separable Banach space, $F$ a bounded subset of its dual $E^*$ and $f \in \overline{\co}^{w^*}(F)$. Then there is
a sequence $(f_n) \subseteq F$ such that
\[\frac{f_1+\dots+f_N}{N} \overset{w^*}{\underset{N \to \infty}{\longrightarrow}} f.\]
\end{prop}

\begin{proof} 
The proof is similar to the proof of Theorem 2(2). We set $X=B_{E^*}$ and assume without loss of generality that 
$F \subseteq X$. Note that, since $E$ is separable, $X$ is weak$^*$ compact and metrizable and also that each
$x \in E$ can be identified with a continuous function on $X$ through the linear isometry
\[T:x \in E \mapsto T(x)=x|X \in C(X) \subseteq \ell_\infty(X).\]
Since the weak$^*$ closed convex hull $\overline{\co}^{w^*}(F) \subseteq X$ is a weak$^*$ compact and metrizable set,
given any $f \in \overline{\co}^{w^*}(F)$ there is a sequence $(\mu_j)$ of convex combinations of elements of $F$
such that 
\begin{equation}
\mu_j \overset{w^*}{\longrightarrow} f.
\end{equation}
We cosider each $\mu_j$ as a finitely supported probability measure on $F \subseteq X$. So if we set $L=T(B_E) \subseteq B_{\ell_\infty(X)}$,
then (7) means that $\mu_j \rightarrow f$ pointwise on $L$.

It then follows from Theorem 1(1) that there is a sequence $(f_n) \subseteq \cup_{j=1}^\infty \supp \mu_j \subseteq F$ such that
\[\frac{\delta_{f_1}+\dots+\delta_{f_N}}{N} \underset{N \to \infty}{\longrightarrow} f \;\; \textrm{pointwise on} \; L;\]
equivalently $\frac{f_1+\dots+f_N}{N} \overset{w^*}{\underset{N \to \infty}{\longrightarrow}} f$.

\end{proof}

\begin{prop}
Let $E$ be a separable Banach space and $\Omega$ be a weak$^*$ compact and convex subset of $E^*$. Then for every $f \in \Omega$
there is a sequence $(f_n) \subseteq \ex \Omega$ such that
\[\frac{f_1+\dots+f_N}{N} \overset{w^*}{\underset{N \to \infty}{\longrightarrow}} f.\]
\end{prop}

\begin{proof} It follows immediately from the Krein-Milman theorem and Proposition 4. 
\end{proof}

Since the dual unit ball $B_{E^*}$ of any Banach space $E$ is a  weak$^*$ compact and convex set, we immediately obtain the following

\begin{cor}
Let $E$ be a separable Banach space. Then for every $f \in B_{E^*}$, there is a sequence $(f_n)$ of extreme points of $B_{E^*}$ such that
\[\frac{f_1+\dots+f_N}{N} \overset{w^*}{\underset{N \to \infty}{\longrightarrow}} f.\]
\end{cor}

An immediate consequence of Proposition 5 is the following well known result.

\begin{theo}
Let $K$ be any compact metric space. Then every Borel probability measure $\mu$ on $K$ (i.e. $\mu \in P(K)$)
admits a u.d. sequence.
\end{theo}

\begin{proof} 
Since $K$ is compact Hausdorff, we have that $P(K)$ is a weak$^*$ compact and convex subset of $C(K)^*=M(K)$
with $\ex P(K)=\{\delta_x:x \in K\}$, thus $P(K)=\overline{\co}^{w^*}(\{\delta_x:x \in K\})$.  $K$ is a metrizable space,
hence $C(K)$ is separable and so Proposition 5 can be applied.

\end{proof}

Applying Corollary 3 to $C(K)^*$, where $K$ is a compact metric space, gives the following result that seems to be new
and will be our motivation for the next section.

\begin{theo}
Let $K$ be any compact metric space. Then for every $\mu \in B=B_{C(K)^*}$ there are sequences $(x_n) \subseteq K$
and $(\varepsilon_n) \subseteq \{\pm 1\}$ such that
\[\frac{\varepsilon_1 \delta_{x_1}+\dots+\varepsilon_N \delta_{x_N}}{N} \overset{w^*}{\underset{N \to \infty}{\longrightarrow}} \mu.\]
(In case when $C(K)$ is the space of continuous complex functions, $(\varepsilon_n) \subseteq \mathbb{T}=\{z \in \mathbb{C}:|z|=1\}).$
\end{theo}

\begin{proof} 
We first assume that $C(K)$ is the space of continuous real functions. We then have
$\ex B=\{\pm \delta_x:x \in K\}$ and hence $B=\overline{\co}^{w^*}(\{\pm \delta_x:x \in K\})$.\\
In the complex case we have that
$\ex B=\{\alpha \delta_x:x \in K \textrm{ and } \alpha \in \mathbb{C},\; |\alpha|=1\}$ and hence $B=\overline{\co}^{w^*}(\{\alpha \delta_x:x \in K 
\textrm{ and } \alpha \in \mathbb{C},\; |\alpha|=1\}\})$ (see Theorem 1.9 of \cite{B}). So the result follows immediately from Corollary 3.

\end{proof}

\begin{rem}
(a) Theorem 3 is of course a direct consequence of Niederreiter's main result (see Theorem 2 of \cite{N}). We state Theorem 3 here, because
it shows that Proposition 5 can be considered as a generalization of such a well known result to the much wider class of separable Banach spaces.

(b) A compact Hausdorff space $K$ is said to be angelic iff for every $A \subseteq K$ and each $x \in \overline{A}$ there is a sequence
$(x_n) \subseteq A$ such that $x_n \rightarrow x$. It is clear that, if the dual unit ball $(B_{E^*},w^*)$ of a Banach space $E$ is an
angelic space, then Proposition 4 and its consequences (Proposition 5, Corollary 3 and Theorems 3,4) remain valid. Well known classes of (not necessarily separable)
Banach spaces with angelic dual balls are weakly compactly generated (WCG) and their generalizations, like weakly countably determined (WCD)
Banach spaces, etc. (see \cite{AM}, \cite{FHHMZ} and \cite{HMVZ}).

We also note that if $K$ is any compact Hausdorff space so that the convex hull $\co(\{\delta_x:x \in K\})$ is weak$^*$ sequentially
dense in $P(K)$, then it is easy to see that every $\mu \in M(K)$ with $\|\mu\|=1$ satisfies the conclusion of Theorem 4 (in particular,
by a result of Niederreiter mentioned in (a), every $\mu \in P(K)$ admits a u.d. sequence).

Finally, assuming Martin's axiom plus the negation of Continuum Hypothesis (MA+$\neg$ CH), Theorem 4 remains valid for every compact separable
space $K$ of topological weight $w(K)<c$, where $c$ is the cardinality of the continuum (the proof is essentially the same as the proof
of Proposition 2.21 of \cite{M2}).

(c) Let $E$ be a Banach space not containing an isomorphic copy of $\ell_1$. Then by a result of Haydon, every weak$^*$ compact and
convex subset $\Omega$ of $E^*$ is the norm closed convex hull of its extreme points (see \cite{H}). It then follows from this result
and the aforementioned considerations that for every $f \in \Omega$ there exists a sequence $(f_n) \subseteq \ex \Omega$ norm
Ces\`{a}ro summable to $f$. It follows in particular that if $E$ is of the form $C(K)$, where $K$ is compact and Hausdorff (since
$\ell_1 \nsubseteq C(K)$, every $\mu \in P(K)$ is purely atomic), then for every $\mu \in P(K)$ there is a $\mu$-u.d. sequence $(x_n)$
in $K$ with the stronger property
\[\| \cdot\|_1-\lim_{N \to \infty} \frac{\delta_{x_1}+\dots+\delta_{x_N}}{N}= \mu.\]
Note that the last result can also be proved by a direct method. We also note that a class of (not necessarily separable)
Banach spaces not containing $\ell_1$ is that of Asplund spaces; a Banach space $E$ is called Asplund, if every separable 
subspace of $E$ has separable dual (see \cite{FHHMZ}).

\end{rem}

\section{Uniformly Distributed sequences with respect to signed measures}

In this section we shall concentrate on duals $C(K)^*=M(K)$ of Banach spaces of the form $C(K)$ with $K$ compact (metrizable) space
and shall further investigate the effect of the results of the previous section to the uniform distribution of sequences in $K$.
We emphasize that everything in this section is over $\mathbb{R}$.

Theorem 4 inspires the following generalization of the classical concept of uniformly distributed sequences defined for
regular Borel probability measures on compact spaces (see Definition 1.1 of \cite{KN}).

\begin{defn}
Let $K$ be a compact Hausdorff space and $\mu \in M(K)$ with total variation $|\mu|(K)=1$ (i.e. $\mu$
is a regular Borel signed measure with $\|\mu\|=1$). We say that a sequence $(x_n) \subseteq K$ is $\mu$-u.d.
if both of the following conditions are satisfied
\begin{enumerate}
\item $(x_n)$ is $|\mu|$-u.d. (in the classical sense) and
\item there is a sequence of signs $(\varepsilon_n) \subseteq \{\pm 1\}$ such that
\[\frac{\varepsilon_1 \delta_{x_1}+\dots+\varepsilon_N \delta_{x_N}}{N} \overset{w^*}{\underset{N \to \infty}{\longrightarrow}} \mu\]
(that is, $ \lim_{N \to \infty} \frac{\varepsilon_1 f(x_1)+\dots+\varepsilon_N f(x_N)}{N} =\int_K f d \mu$, for all $f \in C(K)$).
\end{enumerate}
\end{defn}

We note that:

(a) By applying the above equality for the constant function $f=1$, we get that $\mu(K)=\int_K d \mu=
\lim_{N \to \infty} \frac{\varepsilon_1 +\dots+\varepsilon_N}{N}$ and

(b) if $\mu \in P(K)$ admits a u.d. sequence $(x_n)$, then we can take $\varepsilon_n=1$ for all $n \ge 1$.

Let $\mu \in M(K)$ with $\|\mu\|=1\,(=|\mu|(K))$.  It then follows from Radon-Nikodym theorem that there is a
Borel function $h:K \rightarrow \mathbb{R}$ with $|h(x)|=1$, for all $x \in K$ such that
\[d \mu=h d|\mu|.\]
The function $h$ denoted by $\frac{d \mu}{d|\mu|}$ is the so-called \textsl{Radon-Nikodym derivative} of $\mu$ with respect to its total
variation $|\mu|$. So we have that
\[\int_Kf d \mu=\int_K fh d|\mu|,\]
for all bounded Borel measurable functions $f:K \rightarrow \mathbb{R}$.

With the above notation we have the following.

\begin{prop}
Assume that $|\mu|$ admits a u.d. sequence $(x_n)$ (in the classical sense) and also that the function $h$
is $|\mu|$-Riemann integrable. Then the sequence $(x_n)$ is $\mu$-u.d. (in the sense of Definition 1).
\end{prop}

\begin{proof} 
Let $f \in C(K)$; then the function $fh$ is $|\mu|$-Riemann integrable and since $(x_n)$ is $|\mu|$-u.d.,
we get that (see also Fact II after Remark 7)
\[\int_Kf d \mu=\int_K fh d|\mu|=\lim_{N \to \infty} \frac{(fh)(x_1)+\dots+(fh)(x_N)}{N}. \]
So the desired sequence of signs is the sequence $\varepsilon_n=h(x_n), \; n \ge 1$.

\end{proof}

\begin{rem}
Later in this section, we shall present a class of measures $\mu \in P(K)$, where $K=[a,b]$ is a compact
interval of the real line, such that the function  $h= \frac{d \mu}{d|\mu|}$ is $|\mu|$-Riemann integrable.
\end{rem}

In the sequel we are going to show, essentially by the method of proof of Theorem 2.2, p.183 of \cite{KN}, 
that every measure $\mu \in M(K)$ with $\|\mu\|=1$ admits a u.d. sequence.
We first cite some preliminaries. Let $K$ be a compact Hausdorff space; we denote by $K^\infty$ the cartesian
product of countably many copies of $K$. Then $K^\infty$ is  a compact Hausdorff space endowed with the
product topology. If $\mu \in P(K)$, then $\mu$ induces the product measure $\mu^\infty$ in $K^\infty$,
which we may assume to be complete. We also denote by $\mathcal{B}(K)$ the Banach space of bounded Borel functions
on $K$ endowed with supremum norm.

\begin{theo} 
Let $K$ be a compact metric space and $\mu \in M(K)$ with $\|\mu\|=1$; also let $S$ be the set of all sequences
in $K$, which are $\mu$-u.d. considered as a subset of $K^\infty$. Then $|\mu|^\infty(S)=1$.
\end{theo}

\begin{proof}
We consider a countable total subset $L=\{f_n:n \ge 1\}$ of $C(K)$ with $f_1 \equiv 1$. Also let $h= \frac{d \mu}{d|\mu|}$
(=the Radon-Nikodym derivative of $\mu$ with respect to $|\mu|$).
Set $M=L \cup hL=\{f_n:n \ge 1\} \cup \{f_nh:n \ge 1\}$. As the members of the set $M$ are bounded Borel functions
and $|\mu| \in P(K)$, for each $g \in M$ there is a $|\mu|^\infty$-measurable subset $B_g$ of $K^\infty$ with
$|\mu|^\infty(B_g)=1$ such that
\begin{equation}
\lim_{N \to \infty} \frac{1}{N} \sum_{k=1}^N g(x_k)=\int_K g d |\mu|\;\;\; 
\forall (x_1,\dots,x_k, \dots) \in B_g 
\end{equation}
(see Lemma 2.1, p.182 of \cite{KN}).

Set $B=\cap_{g \in M} B_g$; as the set $M$ is countable, we get that $|\mu|^\infty(B)=1$. Let $(x_1,\dots,x_k, \dots) \in B$.
Since the set $L$ is total in $C(K)$, we get that this sequence is $|\mu|$-u.d. in $K$, that is (8) is valid for every
$f \in C(K)$. Note that the operator $T:f \in C(K) \mapsto hf \in \mathcal{B}(K)$ is a linear isometry, thus the set $hL$
is total in the closed subspace $T(C(K))$ of the Banach space $\mathcal{B}(K)$. So we get that equation (8) also holds
for each member of the space $T(C(K))$, that is
\begin{equation}
\lim_{N \to \infty} \frac{1}{N} \sum_{k=1}^N \varepsilon_k f(x_k)=\int_K fh d |\mu|=\int_K f d\mu\;\;\; 
\forall f \in C(K) 
\end{equation}
where $\varepsilon_k=h(x_k), k=1,2,\dots$. So we are done.

\end{proof}

\begin{rem}
Note that assertions 1. and 2. of Definition 1 (for $f \in C(K)$) are equivalent to the following
\begin{equation}
\lim_{N \to \infty} \frac{(1+\varepsilon_1)f(x_1)+\dots+(1+\varepsilon_N)f(x_N)}{2N}=\int_K f d\mu^+ 
\;\;\textrm{and}
\end{equation}

\begin{equation}
\lim_{N \to \infty} \frac{(1-\varepsilon_1)f(x_1)+\dots+(1-\varepsilon_N)f(x_N)}{2N}=\int_K f d\mu^- 
\end{equation}
where $\mu^+=\frac{1}{2}(|\mu|+\mu)$ and $\mu^-=\frac{1}{2}(|\mu|-\mu)$ are the positive and negative variations
of the measure $\mu$.

Indeed, assuming that 1. and 2. of Definition 1 are valid, for $f \in C(K)$ we have
\[\int_K f d\mu^+=\frac{1}{2} \left[ \int_K f d(|\mu|+\mu) \right]=
\frac{1}{2} \left[ \int_K f d|\mu|+ \int_K f d\mu \right]=\]
\[=\lim_{N \to \infty} \frac{(1+\varepsilon_1)f(x_1)+\dots+(1+\varepsilon_N)f(x_N)}{2N}\]
so (10) holds. In a similar way we get equality (11).

In the converse direction, by adding and subtracting (10) and (11) we get 1. and 2. of Definition 1 respectively.

\end{rem}

Concerning Remark 7, it should be noticed that equalities (10) and (11) (and hence 1. and 2. of Definition 1)
are also valid for every $\mu$-Riemann integrable function (i.e. a bounded function $f:K \rightarrow \mathbb{R}$
that is $|\mu|$-Riemann integrable). In order to prove this, we shall use the following well known facts:

\underline{Fact I} Let $\nu \in M^+(K)$ (with $\nu(K)>0$), then a bounded function $f:K \rightarrow \mathbb{R}$ is
$\nu$-Riemann integrable iff for every $\varepsilon>0$ there are $f_1,f_2 \in C(K)$ such that
\[f_1 \le f \le f_2 \;\, \textrm{and} \; \int_K (f_2-f_1) d\nu \le \varepsilon.\]
(see p.90 of \cite{M2})

\underline{Fact II} Let $\nu \in P(K)$ and $(x_n) \subseteq K$ be a $\nu$-u.d. sequence.
Then for  every $\nu$-Riemann integrable function $f:K \rightarrow \mathbb{R}$
\[\lim_{N \to \infty} \frac{1}{N} \sum_{k=1}^N f(x_k)=\int_K f d\nu.\]

The proofs of Facts I and II are essentially contained in the proofs of Theorem 1.1, p.2 and Theorem 1.2, p.175 (see also
exercise 1.12, p.179) of \cite{KN}.

Let us prove, for instance, equality (10). So let $f:K \rightarrow \mathbb{R}$ be any $\mu$-Riemann integrable function
and $\varepsilon>0$. Then by Fact I there are $f_1 \le f \le f_2$ continuous functions, such that
\[0 \le \int_K (f_2-f_1) d|\mu| \le \varepsilon.\]
Since $0 \le \frac{1 \pm h}{2} \le 1 \; $  (where  $h=\frac{d\mu}{d|\mu|}$) we get that
\[F_1:=f_1 \left( \frac{1+h}{2} \right) \le F:=f \left( \frac{1+h}{2} \right) \le F_2:=f_2 \left( \frac{1+h}{2} \right) \; \textrm{and}\]
\[0 \le \int_K(F_2-F_1) d |\mu|=\int_K(f_2-f_1) \left(\frac{1+h}{2}\right) d|\mu| \le \int_K(f_2-f_1) d |\mu| \le \varepsilon.\]
Set $I=\int_K F d|\mu| \left(=\int_K f \left(\frac{1+h}{2}\right) d|\mu|=\int_K f d\mu^+\right)$.
We then have that
\[I-\varepsilon=\int_K F d|\mu|-\varepsilon \le \int_K F_1 d|\mu|= \int_K f_1 d\mu^+\]
\[=\lim_{N \to \infty} \frac{1}{2N} \sum_{k=1}^N (1+\varepsilon_k) f_1(x_k) \le 
\liminf_{N \to \infty} \frac{1}{2N} \sum_{k=1}^N (1+\varepsilon_k) f(x_k) \]
\[\le \limsup_{N \to \infty} \frac{1}{2N} \sum_{k=1}^N (1+\varepsilon_k) f(x_k) \le 
\lim_{N \to \infty} \frac{1}{2N} \sum_{k=1}^N (1+\varepsilon_k) f_2(x_k)\]
\[= \int_K f_2 d\mu^+=\int_K F_2 d|\mu| \le \int_K F d|\mu|+\varepsilon= I+\varepsilon.\]
Since $\varepsilon$ is arbitrarily small, we get (10).

Taking into account the above remarks and Theorem 5, we obtain the following result.

\begin{theo}
Let $K$ be a compact metric space and $\mu \in M(K)$ with $\|\mu\|=1$. Then there is a sequence $(x_n) \subseteq K$
and a sequence of signs $(\varepsilon_n) \subseteq \{\pm 1\}$ (where $\varepsilon_n=h(x_n)$ and $h=\frac{d\mu}{d|\mu|}$)
such that, for every $\mu$-Riemann integrable function $f:K \rightarrow \mathbb{R}$ we have

\bigskip
$(a) \;\;\; \lim_{N \to \infty} \frac{1}{N} \sum_{k=1}^N f(x_k)= \int_K f d|\mu|;$
in particular $(x_n)$ is $|\mu|$-u.d.

\bigskip
$(b) \;\;\; \lim_{N \to \infty} \frac{1}{N} \sum_{k=1}^N \varepsilon_k f(x_k)= \int_K f d\mu,$

\bigskip
$(c) \;\;\; \lim_{N \to \infty} \frac{1}{2N} \sum_{k=1}^N (1+\varepsilon_k) f(x_k)= \int_K f d\mu^+ \;\, \textrm{and}$

\bigskip
$(d) \;\;\; \lim_{N \to \infty} \frac{1}{2N} \sum_{k=1}^N (1-\varepsilon_k) f(x_k)= \int_K f d\mu^-.$

\end{theo}

In the sequel, we focus on the special case of Theorem 6 when $K$ is a compact interval of the real line, say $K=[a,b]$.
Let $\varphi:[a,b] \rightarrow \mathbb{R}$ be a function of bounded variation. We denote by $\upsilon, p$ and $n$ the (increasing)
functions of total, positive and negative variation of $\varphi\;$ $(\upsilon(x)=V_a^x \varphi, x \in [a,b])$.

We note that these functions are connected as follows:
\[p(x)=\frac{1}{2}(\upsilon(x)+\varphi(x)-\varphi(a)) \;\, \textrm{and}\]
\[n(x)=\frac{1}{2}(\upsilon(x)-\varphi(x)+\varphi(a)),\]
(see p.208 of \cite{Ca}).

We recall that the space $M(K)$ of signed Borel measures on $K$ is in one-to-one correspondence with the space
of functions of bounded variation on $K$ which are right continuous on $(a,b)$ with $\varphi(a)=0$, in the
sense that each $\mu \in M(K)$ is uniquely defined by such a $\varphi$ by the rule
\[\mu((y,x])=\varphi(x)-\varphi(y),\;\; \textrm{for}\; a \le y<x \le b\]
(see Theorem 3.29 of \cite{Fo} and Theorem 14.26 of \cite{Ca}).

With the above notation and terminology, Theorem 6 yields the following

\begin{theo}
Let $\varphi:[a,b] \rightarrow \mathbb{R}$ be a right continuous function (of bounded variation) with total variation $\upsilon(b)=
V_a^b\varphi=1$ and $\varphi(a)=0$. Then there are sequences $(x_n) \subseteq [a,b]$ 
and $(\varepsilon_n) \subseteq \{\pm 1\}$, such that for every point of continuity $x \in [a,b]$ of $\varphi$ we have

\bigskip
$(a') \;\;\; \lim_{N \to \infty} \frac{1}{N} \sum_{k=1}^N \chi_{[a,x)}(x_k)=\upsilon(x)$

\bigskip
$(b') \;\;\; \lim_{N \to \infty} \frac{1}{N} \sum_{k=1}^N \varepsilon_k \chi_{[a,x)}(x_k)=\varphi(x),$

\bigskip
$(c') \;\;\; \lim_{N \to \infty} \frac{1}{2N} \sum_{k=1}^N (1+\varepsilon_k) \chi_{[a,x)}(x_k)=p(x) \;\, \textrm{and}$

\bigskip
$(d') \;\;\; \lim_{N \to \infty} \frac{1}{2N} \sum_{k=1}^N (1-\varepsilon_k) \chi_{[a,x)}(x_k)=n(x).$

\end{theo}

\begin{proof}
Let $\mu=\mu_\varphi$ be the signed (Lebesgue-Stieljes) measure defined by $\varphi$ on $[a,b]$ by the rule
$\mu((a,x])=\varphi(x)\,(=\mu([a,x]))$ for $x \in (a,b]$. As is well known, the Jordan decomposition and the total
variation of $\mu$ are given by
\[\mu=\mu^+-\mu^-,\;\;\; |\mu|=\mu^+ +\mu^-\]
where $|\mu|=\mu_\upsilon,\; \mu^+=\mu_p$ and $\mu^-=\mu_n$, thus in particular $\|\mu\|=\upsilon(b)=1$
(see Theorem 3.29 and exercises 28, 29, p.107 of \cite{Fo}).

Let $D_\varphi$ be the (countable) set of discontinuity points of $\varphi$; then every interval $I \subseteq [a,b]$
whose both endpoints do not belong to $D_\varphi$ has characteristic function which is $\mu$-Riemann integrable. Then
by applying Theorem 6 to intervals of the form $[a,x)$, with $x \notin D_\varphi$ we get the conclusion  (of course 
$\varepsilon_n=h(x_n),\; n \ge 1$ where $h=\frac{d\mu}{d|\mu|}$).

\end{proof}

The last theorem partially generalizes an important result from \cite{KN} (Theorem 4.3, p.138) stating that:

\begin{theo}
Let $\varphi:I=[0,1] \rightarrow \mathbb{R}$ be an increasing function, with $\varphi(0)=0$ and $\varphi(1)=1$.
Then there is a sequence $(x_n) \subseteq I$, such that
\[\lim_{N \to \infty} \frac{1}{N} \sum_{k=1}^N \chi_{[a,x)}(x_k)=\varphi(x),\;\; \textrm{for} \;\, 0 \le x \le 1.\]
We then say that $(x_n)$ has $\varphi$ as the \textsl{asymptotic distribution function mod 1} (abbreviated a.d.f.(mod 1) $\varphi(x)$).
\end{theo}

The proof of this result is given in two steps. First, the continuous case is proved (Lemma 4.2, p.137 of \cite{KN})
and then the general case follows, using a result from real analysis (Lemma 4.3, p.138 of \cite{KN}). We state both of these
results for the reader's convenience.

\begin{lemm}
Let $\varphi:I \rightarrow \mathbb{R}$ be a continuous increasing function, with $\varphi(0)=0$ and $\varphi(1)=1$.
Then there is a sequence $(x_n) \subseteq I$, such that
\[\left| \frac{1}{N} \sum_{k=1}^N \chi_{[0,x)}(x_k)-\varphi(x)\right| \le \frac{\log(N+1)}{N \log 2} \]
for all $N \ge 1$ and $0 \le x \le 1$.
\end{lemm}

Before we state the real analysis lemma, we recall that a continuous function $\varphi:[a,b] \rightarrow \mathbb{R}$,
($a,b \in \mathbb{R}, \; a<b$) is said to be \textsl{polygonal} (or piecewise linear), if its graph consists
of finitely many straight line segments.

\begin{lemm}
Let $\varphi:[a,b] \rightarrow \mathbb{R}$ be an increasing function. Then there is a sequence $(\varphi_k)$ of polygonal increasing
functions defined on $[a,b]$, satisfying $\varphi_k(a)=\varphi(a)$ and $\varphi_k(b)=\varphi(b)$ for $k \ge 1$, which 
converges pointwise to $\varphi$, that is, $\lim_{k \to \infty} \varphi_k(x)=\varphi(x)$, for all $x \in [a,b]$.
\end{lemm}

Our aim is to give a full generalization of both Theorems 7 and 8, in the sense that assertions $(a')$ to $(d')$ of Theorem 7
are valid for every function $\varphi:[a,b] \rightarrow \mathbb{R}$  of bounded variation (with $\varphi(a)=0$ and 
$V_a^b\varphi=1$) and for each point $x \in [a,b]$. We start by generalizing Lemma 4.

\begin{lemm}
Let $\varphi:[a,b] \rightarrow \mathbb{R}$ be a function of bounded variation (with $V_a^b \varphi>0$). Also, let 
$\upsilon, p$ and $n$ be the functions of total, positive and negative variation of $\varphi$. Then there are
sequences of increasing polygonal functions $(g_n)$ and $(h_n)$ defined on $[a,b]$, such that if we let
$\varphi_k=g_k-h_k$, for $k \ge 1$, then we have that ($\varphi_k$ is polygonal and)\\
(i)$\;g_k \rightarrow p,\; h_k \rightarrow n$ and (thus) $\varphi_k \rightarrow \varphi$ pointwise on $[a,b]$;
moreover $g_k(a)=p(a),\; g_k(b)=p(b)$ and $n_k(a)=n(a),\; n_k(b)=n(b)$  for $k \ge 1$.\\
(ii) If $\upsilon_k$ denotes the function of total variation of $\varphi_k$, then $\upsilon_k$ is polygonal
and $\upsilon_k \rightarrow \upsilon$ pointwise on $[a,b]$.
\end{lemm}

\begin{proof}
For each $k \ge 1$ we choose a partition $\mathcal{P}_k=\{t_0^k=a<t_1^k<\dots<t_{m_k}^k=b\}$ of $[a,b]$
with $t_{i+1}^k-t_i^k<\frac{1}{k}$  for $0 \le i <m_k$ that contains all points $x \in (a,b)$ with
$\upsilon(x+0)-\upsilon(x-0)> \frac{1}{k}$ (since $\upsilon$ is increasing, there can only be finitely
many such $x$). We define the functions $g_k$ and $h_k$ as follows: Set $g_k(t_i^k)=p(t_i^k)$ and
$h_k(t_i^k)=n(t_i^k)$  for $0 \le i \le m_k$ and then extend $g_k$ and $h_k$ on $[a,b]$ so as to be linear
on the intervals $[t_i^k,t_{i+1}^k],\; 0 \le i <m_k$. Then clearly $g_k$ and $h_k$ are polygonal
and increasing on $[a,b]$ and (hence) $\varphi_k$ is polygonal on $[a,b]$. 

(i) We shall prove that $(g_k)$ converges pointwise to $p$ (the proof for $(h_k)$ is analogous).
This is trivial for the endpoints $a$ and $b$. Let $x \in (a,b)$; assume first that $x$ is a discontinuity
point of $\upsilon$. Then $\upsilon(x+0)-\upsilon(x-0)>0$ and so, from some $k$ on, we will have that
$x=t_{i_k}^k$ with $0<i_k<m_k$. Therefore $g_k(x)=p(x)$ for sufficiently large $k$.

Now let $\upsilon$ be continuous at $x$, then $\varphi,\;p,\;n$ are also continuous at $x$. So let
$\varepsilon>0$ be given. Then, for all sufficiently large $k$ (say $k \ge k_0$) we will have
\[y \in \left( x-\frac{1}{k},x+\frac{1}{k} \right) \Rightarrow p(y) \in (p(x)-\varepsilon,p(x)+\varepsilon).\]
Yet for each $k$ we have $t_i^k \le x \le t_{i+1}^k$, for some $i=i(k)$ with $0 \le i<m_k$. Since
$0<t_{i+1}^k-t_i^k<\frac{1}{k}$, both $t_i^k,\; t_{i+1}^k$ lie in $ \left( x-\frac{1}{k},x+\frac{1}{k} \right)$.
Hence, for $ k \ge k_0$ we obtain
\[g_k(t_i^k)=p(t_i^k)>p(x)-\varepsilon\;\, \textrm{and}\;\,g_k(t_{i+1}^k)=p(t_{i+1}^k)<p(x)+\varepsilon.\]
Since $g_k$ is increasing, we get that 
\[g_k(t_i^k) \le g_k(x) \le g_k(t_{i+1}^k)\]
and so $p(x)-\varepsilon<g_k(x)<p(x)+\varepsilon$ for $k \ge k_0$, which shows that $g_k(x) \rightarrow p(x)$.

(ii) Let $x \in (a,b]$ (clearly $\upsilon(a)=\upsilon_k(a)=0$ for $k \ge 1$); since $\varphi_k$ is a polygonal
and hence piecewise $C^1$ function, we get that
\[\upsilon_k(x)=\int_a^x |\varphi'_k(t)|dt=\sum_{\lambda=0}^{i(k)-1} |\varphi_k(t_{\lambda+1}^k)-\varphi_k(t_{\lambda}^k)|+
|\varphi_k(t_{i(k)}^k)-\varphi_k(x)|=\]
(where $x \in (t_{i(k)}^k,t_{i(k)+1}^k]$ and $0 \le i(k) <m_k$)

\begin{equation}
=\sum_{\lambda=0}^{i(k)-1} |\varphi(t_{\lambda+1}^k)-\varphi(t_{\lambda}^k)|+
|\varphi(t_{i(k)}^k)-\varphi_k(x)| \le V_a^{t_{i(k)}}\varphi+|\varphi(t_{i(k)}^k)-\varphi_k(x)|.
\end{equation}

It is clear that
\begin{equation}
V_a^{t_{i(k)}}\varphi+|\varphi(t_{i(k)}^k)-\varphi(x)| \le V_a^x\varphi
\end{equation}
since $\varphi_k(x) \rightarrow \varphi(x)$, we get from (12) and (13) that
\begin{equation*}
\limsup_{k \to \infty} \upsilon_k(x) \le \limsup_{k \to \infty} (V_a^{t_{i(k)}}\varphi+|\varphi(t_{i(k)}^k)-\varphi_k(x)|)=
\end{equation*}
\begin{equation}
=\limsup_{k \to \infty} (V_a^{t_{i(k)}}\varphi+|\varphi(t_{i(k)}^k)-\varphi_(x)|) \le V_a^x \varphi.
\end{equation}

But since $\varphi_k \rightarrow \varphi$ pointwise on $[a,b]$, we have that
\begin{equation}
\upsilon(x)=V_a^x \varphi \le \liminf_{k \to \infty} \upsilon_k(x)
\end{equation}
(see exercise 12, p.205 of \cite{Ca}).

It then follows from (14) and (15) that
\[\upsilon(x)=\lim_{k \to \infty} \upsilon_k(x) \;\, \textrm{for}\;\, x \in (a,b].\]
We finally note that it is easy to verify that the function of total variation of a polygonal function is also polygonal.
Therefore each $\upsilon_k$ is a polygonal (and increasing) function.

\end{proof}

We note that the proof of claim (i) of Lemma 5 is similar to the proof of Lemma 4.

\begin{rem}

(1) Regarding the previous Lemma, we set
\[p_k=\frac{1}{2} (\upsilon_k+\varphi_k-\varphi_k(a)) \;\, \textrm{and} \;\,n_k=\frac{1}{2} (\upsilon_k-\varphi_k+\varphi_k(a)), \]
where $\upsilon_k$ is the function of total variation of $\varphi_k$. Then we have that:\\
(i)$\; p_k$ and $n_k$ are the positive and negative variations of $\varphi_k$,\\
(ii) the function $\varphi_k$ is polygonal, the functions $\upsilon_k,\; p_k,\; n_k$ are polygonal and increasing and\\
(iii) $\; \varphi_k \rightarrow \varphi,\; \upsilon_k \rightarrow \upsilon$ pointwise on $[a,b]$ and hence 
$p_k \rightarrow p,\; n_k \rightarrow n$ pointwise on $[a,b]$.

(2) Assume now that $\varphi(a)=0$ and $V_a^b \varphi=1$. We then have that $\varphi_k(a)=0$ for $k \ge 1$ and
$V_a^b \varphi_k=\upsilon_k(b) \underset{k \to \infty}{\longrightarrow} \upsilon(b)=V_a^b \varphi=1$.
Now we define
\[\Phi_k=\frac{\varphi_k}{\upsilon_k(b)},\;\Upsilon_k=\frac{\upsilon_k}{\upsilon_k(b)},\;P_k=\frac{p_k}{\upsilon_k(b)},\;N_k=\frac{n_k}{\upsilon_k(b)} \]
and notice the following:\\
(I) $\; \Upsilon_k,\;P_k,\;N_k$ are the functions of total, positive and negative variation of $\Phi_k$, so that $\Phi_k(a)=0$
and $\Upsilon_k(b)=V_a^b\Phi_k=1$  for $ k \ge 1$.\\
(II) $\; \Phi_k$ is polygonal and $\Upsilon_k,\;P_k,\;N_k$ are polygonal and increasing.\\
(III) For every $x \in [a,b]$ we have that $\Phi_k(x) \rightarrow \varphi(x)$, $\Upsilon_k(x) \rightarrow \upsilon(x)$ and hence
$P_k(x) \rightarrow p(x)$ and $N_k(x) \rightarrow n(x)$.

\end{rem}

It follows from the aforementioned remark that Lemma 5 can be stated as follows:

\begin{prop}
Let $\varphi:[a,b] \rightarrow \mathbb{R}$ be a function of bounded variation with $\varphi(a)=0$ and  $V_a^b \varphi=1$.
Also, let $\upsilon, p$ and $n$ be the functions of total, positive and negative variation of $\varphi$. Then there is 
a sequence $\varphi_k:[a,b] \rightarrow \mathbb{R}, \; k \ge 1$ of polygonal functions, such that if $\upsilon_k, p_k$ and $n_k$ 
are the total, positive and negative variations of $\varphi_k$, then (these functions are polygonal and)\\
(1) $\;\varphi_k \rightarrow \varphi$, $\upsilon_k \rightarrow \upsilon$, $p_k \rightarrow p$ and $n_k \rightarrow n$
pointwise on $[a,b]$.\\
(2) $\; \varphi_k(a)=0$ and $V_a^b \varphi_k=1$  for $k \ge 1$.

\end{prop}

Let $\varphi:I=[0,1] \rightarrow \mathbb{R}$ be a continuous function of bounded variation with $\varphi(0)=0$ 
and  $V_0^1 \varphi=1$. Denote, as usual, the function of total variation of $\varphi$ by $\upsilon$. Then
by Theorem 7, there are sequences $\omega=(x_n) \subseteq I$ and $\varepsilon=(\varepsilon_n) \subseteq \{\pm1\}$
such that\\
\bigskip
$(i) \;\;\; \lim_{N \to \infty} \frac{1}{N} \sum_{k=1}^N \chi_{[0,x)}(x_k)=\upsilon(x)$ and\\
\bigskip
$(ii) \;\;\; \lim_{N \to \infty} \frac{1}{N} \sum_{k=1}^N \varepsilon_k \chi_{[0,x)}(x_k)=\varphi(x) $\\
for all $0 \le x \le 1$.

We can now define the \textsl{discrepancy} $D_N(\omega;\upsilon)$ of $\omega=(x_n)$ with respect to the (continuous) function $\upsilon$
by the rule
\[D_N(\omega;\upsilon)=\sup_{0 \le a<b \le 1} \left|\frac{1}{N} \sum_{k=1}^N \chi_{[a,b)} (x_k)-(\upsilon(b)-\upsilon(a))\right|,\]
and one can prove that $\lim_{N \to \infty} D_N(\omega;\upsilon)=0$ (see Theorem 1.1, p.89 and the remarks before Theorem 1.2, p.90 of \cite{KN}).

Similarly, we define the discrepancy of $\omega$ with respect to $\varphi$ as follows (recall that $(\varepsilon_k) \subseteq \{\pm 1\}$)
\[D_N(\omega;\varphi)=\sup_{0 \le a<b \le 1} \left|\frac{1}{N}  \sum_{k=1}^N \varepsilon_k \chi_{[a,b)} (x_k)-(\varphi(b)-\varphi(a))\right|.\]
Then an analogous result can be shown, the proof of which is similar to the proof of Theorem 1.1, p.89 of \cite{KN}. So, with the above assumptions 
and notation for $\varphi$ we have the following

\begin{lemm}
$\lim_{N \to \infty} D_N(\omega;\varphi)=0$ 
\end{lemm}

\begin{proof} 
We first define the discrepancies of the functions $p$ and $n$ in the obvious way. It is easy to see that
\begin{equation}
D_N(\omega;\varphi) \le D_N(\omega;p)+D_N(\omega,n).
\end{equation}
Since $p$ and $n$ are continuous and hence uniformly continuous on the compact interval $I$, given any $\varepsilon>0$
there is $m \in \mathbb{N}$ such that
\begin{equation}
x,y \in I \;\, \textrm{and} \;\, |x-y|<\frac{1}{m} \Rightarrow |p(x)-p(y)|<\frac{\varepsilon}{2} \;\, \textrm{and}\;\, |n(x)-n(y)|<\frac{\varepsilon}{2}.
\end{equation}
We may pick $m$ so large that $\frac{1}{m}<\varepsilon$.

For such an integer $m$, set $I_k=\left[ \frac{k}{m}, \frac{k+1}{m} \right),\; 0 \le k \le m-1$. Using the equalities $(c')$ and $(d')$ of Theorem 7, we get an $N_0=N_0(m) \in \mathbb{N}$ such that for every $N \ge N_0$ and each $k=0,1,\dots,m$ we have
\begin{equation}
\mu^+(I_k)-\frac{1}{m^2} \le \frac{1}{N} \sum_{\lambda=1}^N \left(\frac{1+\varepsilon_\lambda}{2}\right) \chi_{I_k}(x_\lambda) \le \mu ^+(I_k)+\frac{1}{m^2}
\end{equation}

\begin{equation}
\textrm{and} \;\,\mu^-(I_k)-\frac{1}{m^2} \le \frac{1}{N} \sum_{\lambda=1}^N \left(\frac{1-\varepsilon_\lambda}{2}\right) \chi_{I_k}(x_\lambda) \le \mu ^-(I_k)+\frac{1}{m^2}
\end{equation}
where $\mu^+(I_k)=p \left( \frac{k+1}{m} \right)-p\left( \frac{k}{m} \right)$, $\mu^-(I_k)=n \left( \frac{k+1}{m} \right)-n\left( \frac{k}{m} \right)$
and $\mu^+=\mu_p$, $\mu^-=\mu_n$ are the positive and negative variations of $\mu=\mu_\varphi$ (cf. the proof of Theorem 7).

Now consider an arbitrary interval $J=[a,b] \subseteq I$, then there are subintervals $J_1,\; J_2$ of $J$ each one of them being a finite union of 
succesive intervals $I_k$, such that $J_1 \subseteq J \subseteq J_2$ and
\begin{equation}
\mu^+(J)-\mu^+(J_1)<\varepsilon,\; \mu^+(J_2)-\mu^+(J)< \varepsilon \;\, \textrm{and}\;\,\mu^-(J)-\mu^-(J_1)<\varepsilon,\; 
\mu^-(J_2)-\mu^-(J)< \varepsilon.
\end{equation}

These inequalities are easy consequences of (17), i.e. of the uniform continuity of $p$ and $n$.
By adding at most $m$ inequalities of the form (18), we get that
\[\mu^+(J_1)-\frac{1}{m} \le \frac{1}{N} \sum_{\lambda=1}^N \left(\frac{1+\varepsilon_\lambda}{2}\right) \chi_{J_1}(x_\lambda) \le
\frac{1}{N} \sum_{\lambda=1}^N \left(\frac{1+\varepsilon_\lambda}{2}\right) \chi_J(x_\lambda) \le \]
\[\le \frac{1}{N} \sum_{\lambda=1}^N  \left(\frac{1+\varepsilon_\lambda}{2}\right) \chi_{J_2}(x_\lambda) \le \mu^+(J_2)+\frac{1}{m};\]
then using (20) we conclude that
\begin{equation}
\mu^+(J)-2 \varepsilon<\mu^+(J)-\,\frac{1}{m}-\varepsilon \le  \frac{1}{N} \sum_{\lambda=1}^N \left(\frac{1+\varepsilon_\lambda}{2}\right) \chi_J(x_\lambda) \le
\mu^+(J)+\,\frac{1}{m}+\varepsilon < \mu^+(J)+2 \varepsilon. 
\end{equation}
In a similar way we get that
\begin{equation}
\mu^-(J)-2 \varepsilon<\mu^-(J)-\,\frac{1}{m}-\varepsilon \le  \frac{1}{N} \sum_{\lambda=1}^N \left(\frac{1-\varepsilon_\lambda}{2}\right) \chi_J(x_\lambda) \le
\mu^-(J)+\,\frac{1}{m}+\varepsilon < \mu^-(J)+2 \varepsilon. 
\end{equation}
Since (21) and (22) are independent of $J$, we conclude that
$\lim_{N \to \infty} D_N(\omega;p)=0, \;\;$ $\lim_{N \to \infty} D_N(\omega;n)=0$ and thus by (16)  $\lim_{N \to \infty} D_N(\omega;\varphi)=0$.

\end{proof}

\begin{rem}
It is also possible (and useful) to define the concept of $D_N^*$ discrepancy for the functions $\varphi,\; \upsilon,\;p$ and $n$ 
(cf. Definition 1.2, p.90 of \cite{KN}). For instance we may define
\[D_N^*(\omega;\varphi)=\sup_{0 \le x \le 1} \left|\frac{1}{N}  \sum_{k=1}^N \varepsilon_k \chi_{[0,x)} (x_k)-\varphi(x) \right|.\]
It is then easy to see that
\[D_N^*(\omega;\varphi) \le D_N(\omega;\varphi) \le 2D_N^*(\omega;\varphi) \]
(cf. Theorem 1.3, p.91 of \cite{KN}). So we get that
\[\lim_{N \to \infty} D_N(\omega;\varphi)=0 \Leftrightarrow \lim_{N \to \infty} D_N^*(\omega;\varphi)=0. \]
\end{rem}

Let $(f_n)$ be a sequence of scalar valued functions defined on a set $X$. Given a strictly increasing sequence of positive integers
$1 \le N_1<N_2<\dots<N_k<\dots$, we arrange the terms of $(f_n)$ setting
\[g_N=f_1  \;\, \textrm{for} \;\, 1 \le N < N_1 \;\, \textrm{and} \;\, g_N=f_k, \;\, \textrm{for} \;\, N_{k-1} \le N < N_k,\; k \ge 2.\]
Then the following lemma has an easy proof, which we omit.

\begin{lemm}
If $f_n \rightarrow f$ pointwise on $X$, then $g_N \rightarrow f$ pointwise on $X$.
\end{lemm}

Let now $f_n:\mathbb{N} \rightarrow \mathbb{R}$, $n \ge 1$ be a sequence of functions, such that for every $n \ge 1$, $\lim_{m \to \infty} f_n(m)=0$.
We consider a strictly increasing sequence of positive integers $(N_k)_{k \ge 1}$ such that
\[m \ge N_k \Rightarrow |f_k(m)| \le \frac{1}{k}  \;\, \textrm{for}\;\, k \ge 1.\]
(Since each $f_n$ is a null sequence of scalars, such a sequence exists). If we apply the above arrangement to $(f_n)$ defined by $(N_k)$, 
we get the following.

\begin{lemm}
$\lim_{N \to \infty} g_N(N)=0$.
\end{lemm}

\begin{proof} 
Let $N \ge N_1$, then there is $k \ge 2: N_{k-1} \le N < N_k$, hence $g_N(N)=f_k(N)$. But since $N \ge N_k$, we get (from the definition
of $(N_k)$) that $|f_k(N)| \le \frac{1}{k}$ and so
\[|g_N(N)|=|f_k(N)| \le \frac{1}{k}.\]
As $N \rightarrow \infty$ implies $k \rightarrow \infty$, we obtain the desired result.

\end{proof}

We are now in a position to prove the desired generalization of Theorems 7 and 8.

\begin{theo}
Let $\varphi:[a,b] \rightarrow \mathbb{R}$ be a function of bounded variation with $\varphi(a)=0$ and $V_a^b \varphi=1$.
Also let $\upsilon,\;p$ and $n$ be the functions of total, positive and negative variations of $\varphi$.
Then there are sequences $\tau=(x_n) \subseteq [a,b]$ 
and $\varepsilon=(\varepsilon_n) \subseteq \{\pm 1\}$, such that for every $x \in [a,b]$ we have

\bigskip
$(a'') \;\;\; \lim_{N \to \infty} \frac{1}{N} \sum_{k=1}^N \chi_{[a,x)}(x_k)=\upsilon(x),$

\bigskip
$(b'') \;\;\; \lim_{N \to \infty} \frac{1}{N} \sum_{k=1}^N \varepsilon_k \chi_{[a,x)}(x_k)=\varphi(x),$

\bigskip
$(c'') \;\;\; \lim_{N \to \infty} \frac{1}{2N} \sum_{k=1}^N (1+\varepsilon_k) \chi_{[a,x)}(x_k)=p(x), \;\, \textrm{and}$

\bigskip
$(d'') \;\;\; \lim_{N \to \infty} \frac{1}{2N} \sum_{k=1}^N (1-\varepsilon_k) \chi_{[a,x)}(x_k)=n(x).$

\end{theo}

\begin{proof} 
We first reduce the theorem to the case when $[a,b]$ is the unit interval $I=[0,1]$. Consider the affine continuous function
$g(x)=(b-a)x +a, \; x \in I$; clearly $g$ is strictly increasing, with $g(0)=a$ and $g(1)=b$. Set $\Phi=\varphi \circ g$ and notice that\\
(i) $\; \Phi$ is of bounded variation on $I$ with $\Phi(0)=0$ and $V_0^1\Phi=V_a^b\varphi=1$ and\\
(ii) $\;\upsilon_\Phi(x)=\upsilon_\varphi(g(x)),\; p_\Phi(x)=p_\varphi(g(x))$ and $n_\Phi(x)=n_\varphi(g(x))$, for $x \in I$. Now let
$(z_n) \subseteq I$ and $(\varepsilon_n) \subseteq \{\pm 1\}$ satisfying conditions $(a'')$ to $(d'')$ for the function $\Phi$. Then the
sequences $x_n=g(z_n), n \ge 1$ and $(\varepsilon_n)$ satisfy the same conditions for $\varphi$.

In order to prove the theorem (with $[a,b]=I$) we will follow the method of proof of Theorem 8 (Theorem 4.3, p.138 of \cite{KN}) and use Proposition 7.
So let $\varphi_k,\;\upsilon_k$ (and $p_k,\; n_k$) be as in Proposition 7. Since each $\upsilon_k$ is continuous and increasing with  $\upsilon_k(0)=0$
and $\upsilon_k(1)=1$, by Lemma 3 (Lemma 4.2, p.137 of \cite{KN}) there is a sequence $\tau_k=(x_1^k,x_2^k,\dots,x_n^k,\dots)$ satisfying
\begin{equation}
\left| \frac{1}{N} \sum_{n=1}^N \chi_{[0,x)}(x_n^k)-\upsilon_k(x) \right| \le \frac{\log(N+1)}{N \cdot \log 2}
\end{equation}
for all $N \ge 1$ and $x \in I$.

It is  clear that if we fix some $k \in \mathbb{N}$, then letting $N \rightarrow \infty$ we have
\begin{equation}
\lim_{N \to \infty} \frac{1}{N} \sum_{n=1}^N \chi_{[0,x)}(x_n^k)=\upsilon_k(x)\;\, \textrm{for} \;\, x \in I.
\end{equation}

Let $\mu_k=\mu_{\varphi_k}$ be the Lebesgue-Stieljes measure that $\varphi_k$ defines on $I$ and $h_k=\frac{d\mu_k}{d|\mu_k|}$
be the corresponding Radon-Nikodym derivative. Since each $\varphi_k$ is polygonal, its derivative is a step function, hence
$h_k$ is $|\mu_k|$-Riemann integrable, which implies by Proposition 6 that the sequence $\varepsilon_k=(\varepsilon_n^k)_{n \ge 1}$,
where $\varepsilon_n^k=h_k(x_n^k), n \ge 1$ has the property
\begin{equation}
\lim_{N \to \infty} \frac{1}{N} \sum_{n=1}^N \varepsilon_n^k \chi_{[0,x)}(x_n^k)=\varphi_k(x)\;\, \textrm{for} \;\, x \in I.
\end{equation}
It follows from Lemma 6 and Remark 9 that, if we set
\[D_N^*(\tau_k;\varphi_k)=\sup_{0 \le x \le 1} \left|\frac{1}{N}  \sum_{n=1}^N \varepsilon_n^k \chi_{[0,x)} (x_n^k)-\varphi_k(x) \right|,\]
then we have 
\begin{equation}
D_N^*(\tau_k;\varphi_k) \underset{N \to \infty}{\longrightarrow}0 \;\, \textrm{for every}\;\, k \ge 1.
\end{equation}
We notice that we may furthermore assume that 
\begin{equation}
D_N^*(\tau_N;\varphi_N) \underset{N \to \infty}{\longrightarrow}0.
\end{equation}
In order to obtain (27), we consider a strictly increasing sequence of positive integers $(N_k)_{k \ge 1}$ such that
\[m \ge N_k \Rightarrow D_m^*(\tau_k;\varphi_k) \le \frac{1}{k} \;\, \textrm{for} \;\, k \ge 1.\]
Then we arrange the sequence of functions $(\varphi_k)_{k \ge 1}$ as in Lemma 7, that is we set
\[\Phi_N=\varphi_1 \;\, \textrm{for} \;\, 1 \le N < N_1 \;\, \textrm{and} \;\, \Phi_N=\varphi_k \;\, \textrm{for} \;\, N_{k-1} \le N < N_k,\; k \ge 2.\]
It then follows from Lemmas 7 and 8 that $\Phi_k \rightarrow \varphi$, ($ \Upsilon_k \rightarrow \upsilon$, etc.) pointwise on $I$ and that
(23) to (27) remain valid for the sequence $(\Phi_k)_{k \ge 1}$. So we may (and will) assume without loss of generality that
$\Phi_k=\varphi_k$  for $k \ge 1$.

Now we construct the sequence $\tau=(x_n) \subseteq I$ by listing succesively the first term of $\tau_1$, the first two terms of $\tau_2$
, $\dots$, the first $k$ terms of $\tau_k$, that is,
\[\tau=(x_1^1,x_1^2,x_2^2,\dots,x_1^k,x_2^k,\dots,x_k^k,\dots).\]
The sequence of signs $\varepsilon=(\varepsilon_n)$ is constructed similarly; so we set
\[\varepsilon=(\varepsilon_1^1,\varepsilon_1^2,\varepsilon_2^2,\dots,\varepsilon_1^k,\varepsilon_2^k,\dots,\varepsilon_k^k,\dots).\]
We are going to prove that $\tau$ and $\varepsilon$ have the desired properties. Assertion $(a'')$ of this theorem is proved the same 
way as in the proof of Theorem 8.

Indeed, by Lemma 4.1, p.136 of \cite{KN}, it suffices to prove that
\[ \lim_{k \to \infty} \frac{1}{k} \sum_{i=1}^k \chi_{[0,x)}(x_i^k)=\upsilon(x) \;\, \textrm{for}\;\, x \in I.\]
By using (23) and as (by Proposition 7) $\upsilon_k(x) \rightarrow \upsilon(x)$  for $x \in I$, we conclude that
\[\left| \frac{1}{k} \sum_{i=1}^k \chi_{[0,x)}(x_i^k)-\upsilon(x) \right| \le \left| \frac{1}{k} \sum_{i=1}^k \chi_{[0,x)}(x_i^k)-\upsilon_k(x) \right|
+|\upsilon_k(x)- \upsilon(x)| \le\]
\[\le \frac{\log(k+1)}{k \log 2}+|\upsilon_k(x)- \upsilon(x)| \underset{k \to \infty}{\longrightarrow} 0 \;\, \textrm{for} \;\, x \in I.\]
This way assertion $(a'')$ is proved.

Now, to prove assertion $(b'')$, using again Lemma 4.1 of \cite{KN} as above it suffices to show that
\[\lim_{k \to \infty} \frac{1}{k} \sum_{i=1}^k \varepsilon_i^k \chi_{[0,x)}(x_i^k)=\varphi(x)  \;\, \textrm{for} \;\, x \in I.\]
We now use (27) and the fact that $\varphi_k(x) \rightarrow \varphi(x)$  for $x \in I$ (see Proposition 7), so we get
\[\left| \frac{1}{k} \sum_{i=1}^k \varepsilon_i^k \chi_{[0,x)}(x_i^k)-\varphi(x) \right| \le \left| \frac{1}{k} 
\sum_{i=1}^k \varepsilon_i^k \chi_{[0,x)}(x_i^k)-\varphi_k(x) \right|
+|\varphi_k(x)- \varphi(x)| \le\]
\[\le D_k^*(\tau_k;\varphi_k) +|\varphi_k(x)- \varphi(x)| \underset{k \to \infty}{\longrightarrow} 0.\]
Assertions $(c'')$ and $(d'')$ follow easily from $(a'')$ and $(b'')$. The proof of the theorem is now complete.

\end{proof}

\noindent \textbf{Concluding remarks}\\
1. Let $\varphi:[a,b] \rightarrow \mathbb{R}$ be a differentiable function with ($\varphi(a)=0$ and) bounded derivative.
Then it is Lipschitz continuous and (hence) of bounded variation. Let $\mu=\mu_\varphi$ be the Lebesgue-Stieljes measure
defined by $\varphi$ on $[a,b]$ and $h=\frac{d \mu}{d|\mu|}$. Assuming that $\varphi'$ is Riemann integrable, it is not 
difficult to show that $h$ is $|\mu|$-Riemann integrable. It then easily follows that if $\varphi$ is piecewise $C^1$
(for instance a polygonal function), then $h$ is $|\mu|$-Riemann integrable.

On the other hand if $\varphi'$ is not Riemann integrable, that is, $\varphi$ is a Volterra type function, then the function $h$
may or may not be $|\mu|$-Riemann integrable. For examples (and the properties) of functions of Volterra type, we refer
the reader to the books \cite{Go}, pp.35-36 and \cite{Br}, pp.22-25 and 33-35.

\bigskip
\noindent 2. Concerning future work, we note the following:\\
(a) It would be interesting to have a generalization of Theorems 7 and 8 for functions of bounded variation of several variables,
that is, for functions $f$ defined on the cube $I^n$ for $n \ge 2$, which satisfy a proper notion of bounded variation
(see for instance Definition 5.2, p.147 of \cite{KN}).

\bigskip
\noindent (b) Besides compact metric spaces (see Theorem 3), there are several classes of compact non-metrizable spaces $K$ with the
property that every measure $\mu \in P(K)$ admits a u.d. sequence (see \cite{Lo} and \cite{M2}). For such spaces it would be interesting
to know if every signed measure $\mu$ with $\|\mu\|=1$ admits a u.d. sequence in the sense of Definition 1 (cf. also Remark 5(b)).
In our opinion the most interesting case is that of compact separable groups $G$; since we know that under special set-theoretic
assumptions, i.e. Continuum Hypothesis (CH) (see \cite{Lo}) or Martin's axiom plus the negation of Continuum Hypothesis (MA+$\neg$ CH), 
(see \cite{FrPl}), every measure $\mu \in P(G)$ admits a u.d. sequence. We note that it is enough to consider the compact group $\{0,1\}^c$,
where $c=$ the cardinality of the continuum; this is so because, as is well known, every compact separable group is a dyadic space,
i.e. a continuous image of $\{0,1\}^c$ (see \cite{GM} and \cite{M2}).

\bigskip
\noindent (c) Regarding our consideration, the case of locally compact and separable metrizable spaces is also of interest, see \cite{KN}
Notes, pp.177-178.

\section*{Declarations}

The authors have no relevant financial or non-financial interests to disclose.

\normalsize

\end{document}